\documentclass[11pt]{article}
\usepackage{}
\usepackage{mathrsfs}
\usepackage{amsfonts}

\usepackage{CJK}
\usepackage{amsfonts,amssymb,amsmath,mathrsfs}

\usepackage{color,latexsym,amsfonts}
 \newtheorem{theorem}{Theorem}[section]
 \newtheorem{lemma}[theorem]{Lemma}
 \newtheorem{corollary}[theorem]{Corollary}
 \newtheorem{proposition}[theorem]{Proposition}
 
 \newtheorem{Definition}[theorem]{Definition}
 \newtheorem{remark}[theorem]{Remark}
 \newtheorem{condition}[theorem]{Condition}

 \def\blemma{\begin{lemma}\sl{}\def\elemma{\end{lemma}}}
 \def\btheorem{\begin{theorem}\sl{}\def\etheorem{\end{theorem}}}
 \def\bcorollary{\begin{corollary}\sl{}\def\ecorollary{\end{corollary}}}
 
 \def\bproposition{\begin{proposition}\sl{}\def\eproposition{\end{proposition}}}

 \def\beqlb{\begin{eqnarray}}\def\eeqlb{\end{eqnarray}}
 \def\beqnn{\begin{eqnarray*}}\def\eeqnn{\end{eqnarray*}}

 \def\<{\langle}\def\>{\rangle}

 \def\ar{&\!\!}

 \def\eqref#1{{\rm(\ref{#1})}}

\def\d{\textup{d}}
\def\I{\textup{I}}
\def\e{\textup{e}}

\def\fin{\hfill$\square$}

\def\newdot{{\kern.8pt\cdot\kern.8pt}}

\def\B{\mathscr{B}}

\def\R{\mathbb{R}}
\def\E{\mathbb{E}}
\def\P{\mathbb{P}}

\def\<{\langle}
\def\>{\rangle}
\def\Proof.{\noindent{\bf Proof.}}

\begin{document}

\

\noindent{}

\bigskip\bigskip

\centerline{\Large\bf Integration by parts formula and applications for SDE}

\smallskip

\centerline{\Large\bf driven by fractional Brownian motion
\footnote{Supported by the Research project of Natural Science Foundation of
   Anhui Provincial Universities (Grant No. KJ2013A134), National Natural Science Foundation of China
(Grant No. 11371029).}
}

\smallskip
\
\bigskip\bigskip

\centerline{Xiliang Fan}

\bigskip

\centerline{Department of Statistics, Anhui Normal University, Wuhu 241003, China}
\smallskip
\centerline{fanxiliang0515@163.com}

\smallskip

\bigskip\bigskip

{\narrower{\narrower

\noindent{\bf Abstract.} By constructing a new family of successful couplings, 
the Driver-type integration by parts formula is established for the operator associated with
stochastic differential equation driven by fractional Brownian motion.
As applications, shift Harnack type inequalities are presented
and then the absolute continuity of the solution is proved.
}}


\bigskip



\section{Introduction}

\setcounter{equation}{0}

The Driver integration by parts formula due to \cite{Driver97} is a powerful tool for the underlying Markov semigroups.
This, together with the Bismut derivative formula \cite{Bismut84}, allows us to derive regular estimates for commutator,
which plays a key role in the study of flow properties \cite{Fang&Li&Luo11}.
Recently, based upon the back coupling method, Wang \cite{Wang12b} established the integration by parts formulae 
for various models including degenerate diffusion processes, delayed SDEs and semi-linear SPDEs,
and showed that, in general the integration by parts formula is more complicated and harder to obtain than the
derivative formula.
Afterwards, 
Zhang \cite{ZhangSQ12} studied semi-linear SPDE with delay; 
Fan \cite{Fan13b} discussed stochastic Volterra equation;
Wang \cite{Wang13a} considered SDE with L\'{e}vy noise.

In this paper, we are concerned with SDEs driven by fractional Brownian motion.
It is well known that this kind of noise is not Markovian and even more not semimartingale.
Now, there exist numerous attempts to define a stochastic integral with respect to the fractional Brownian motion 
and then many works to discuss SDEs with the noise.
For instance,
based on the approach of \cite{Zahle98a}, Nualart and R\u{a}\c{s}canu \cite{Nualart&Rascanu02a} proved the existence and uniqueness result with Hurst parameter $H>1/2$;
Coutin and Qian \cite{Coutin&Qian02a} also derived the existence and uniqueness result for $H\in(1/4,1/2)$ via the theory of rough path analysis introduced in \cite{Lyons98a};
in \cite{Hairer&Ohashi07,Hairer&Pillai11} and \cite{Saussereau12}, the authors studied the ergodicity and Talagrand's transportation inequalities for the solutions, respectively;
Fan \cite{Fan13,Fan14} established Harnack inequalities for the solution with $H<1/2$ and $H>1/2$, respectively.
As for the regularities, the readers may refer to \cite{Baudoin&Ouyang11,Hu&Nualart&Song08,Nourdin&Simon06,Nualart&Saussereau09} and references therein.
However, as far as we know, the study of the existence of density of the solution mainly depends on the Malliavin calculus.
Motivated by the work \cite{Wang12b}, 
we will be able to state the absolute continuity of the law as a consequence of shift Harnack type inequalities,
which will be implied by the integration by parts formula.

The paper is organized as follows.
In the next section we give some preliminaries on fractional Brownian motion.
The integration by parts formula is discussed in section 3.
Then, in section 4 we investigate some applications.

\section{Preliminaries}

\setcounter{equation}{0}

In this part, we will recall some basic results about fractional Brownian motion.
The main references for all these results are \cite{Alos&Mazet&Nualart01a}, \cite{Biagini&Hu08a}, \cite{Decreusefond&Ustunel98a}.

Let $H\in(0,1)$.
The $d$-dimensional fractional Brownian motion with Hurst parameter $H$ on the probability space $(\Omega,\mathscr{F},\mathbb{P})$
can be defined as the centered Gauss process $B^H=\{B_t^H, t\in[0,T]\}$ with covariance function
$\E \left(B_t^{H,i}B_s^{H,j}\right)=R_H(t,s)\delta_{i,j}$,
where
\beqnn
R_H(t,s)=\frac{1}{2}\left(t^{2H}+s^{2H}-|t-s|^{2H}\right).
\eeqnn
In particular, if $H=1/2, B^H$ is a $d$-dimensional Brownian motion.
Furthermore, one can show that $\E|B_t^{H,i}-B_s^{H,i}|^p=C(p)|t-s|^{pH},\ \forall p\geq 1$.
Consequently, $B^{H,i}$ have $(H-\epsilon)$-H\"{o}lder continuous paths for all $\epsilon>0,\ i=1,\cdot\cdot\cdot,d$.

For each $t\in[0,T]$, let $\mathcal {F}_t$ be the $\sigma$-algebra generated by the random variables $\{B_s^H:s\in[0,t]\}$ and the
$\mathbb{P}$-null sets.

Denote $\mathscr{E}$ by the set of step functions defined on $[0,T]$.
Let $\mathcal {H}$ be the Hilbert space defined as the closure of $\mathscr{E}$ with respect to the scalar product
\beqnn
\langle (I_{[0,t_1]},\cdot\cdot\cdot,I_{[0,t_d]}),(I_{[0,s_1]},\cdot\cdot\cdot,I_{[0,s_d]})\rangle_\mathcal {H}=\sum\limits_{i=1}^dR_H(t_i,s_i).
\eeqnn
By B.L.T. theorem, the mapping $(I_{[0,t_1]},\cdot\cdot\cdot,I_{[0,t_d]})\mapsto\sum_{i=1}^dB_{t_i}^{H,i}$ can be extended to an isometry
between $\mathcal {H}$ and the Gauss space $\mathcal {H}_1$ associated with $B^H$.
We denote the isometry between $\mathcal {H}$ and $\mathcal {H}_1$ by $\phi\mapsto B^H(\phi)$.
On the other hand, the covariance kernel $R_H(t,s)$ can be written as
\beqnn
R_H(t,s)=\int_0^{t\wedge s}K_H(t,r)K_H(s,r)\d r,
\eeqnn
where $K_H$ is a square integrable kernel given by
\beqnn
K_H(t,s)=\Gamma\left(H+\frac{1}{2}\right)^{-1}(t-s)^{H-\frac{1}{2}}F\left(H-\frac{1}{2},\frac{1}{2}-H,H+\frac{1}{2},1-\frac{t}{s}\right),
\eeqnn
in which $F(\cdot,\cdot,\cdot,\cdot)$ is the Gauss hypergeometric function (for details see \cite{Nikiforov&Uvarov88}).

Now, we define the linear operator $K_H^*:\mathscr{E}\rightarrow L^2([0,T],\R^d)$ by
\beqnn
(K_H^*\phi)(s)=K_H(T,s)\phi(s)+\int_s^T(\phi(r)-\phi(s))\frac{\partial K_H}{\partial r}(r,s)\d r.
\eeqnn
In \cite{Alos&Mazet&Nualart01a}, the authors showed that for all $\phi,\psi\in\mathscr{E}$,
\beqnn
\langle K_H^*\phi,K_H^*\psi\rangle_{L^2([0,T],\R^d)}=\langle\phi,\psi\rangle_\mathcal {H},
\eeqnn
and therefore $K_H^*$ is  an isometry between $\mathcal{H}$ and $L^2([0,T],\R^d)$.
Consequently, the fractional Brownian motion $B^H$ has the following integral representation
\beqnn
 B^H_t=\int_0^tK_H(t,s)\d W_s,
\eeqnn
where $\{W_t=B^H((K_H^*)^{-1}{\rm I}_{[0,t]})\}_{t\in[0,T]}$ is a Wiener process.

Consider the operator $K_H: L^2([0,T],\R^d)\rightarrow I_{0+}^{H+1/2 }(L^2([0,T],\R^d))$ associated with the integrable kernel $K_H(\cdot,\cdot)$
\beqnn
(K_Hf^i)(t)=\int_0^tK_H(t,s)f^i(s)\d s,\ i=1,\cdot\cdot\cdot,d.
\eeqnn
By \cite{Decreusefond&Ustunel98a}, we know that $K_H$ is an isomorphism and moreover,  for each $f\in L^2([0,T],\mathbb{R}^d)$,
\beqnn
(K_H f)(s)=I_{0+}^{2H}s^{1/2-H}I_{0+}^{1/2-H}s^{H-1/2}f,\ H\leq1/2,
\eeqnn
\beqnn
(K_H f)(s)=I_{0+}^{1}s^{H-1/2}I_{0+}^{H-1/2}s^{1/2-H}f,\ H\geq1/2.
\eeqnn
Therefore, for each $h\in I_{0+}^{H+1/2}(L^2([0,T],\R^d))$, the inverse operator $K_H^{-1}$ is of the form
\beqlb\label{2.0'}
(K_H^{-1}h)(s)=s^{H-1/2}D_{0+}^{H-1/2}s^{1/2-H}h',\ H>1/2,
\eeqlb
\beqlb\label{2.0''}
(K_H^{-1}h)(s)=s^{1/2-H}D_{0+}^{1/2-H}s^{H-1/2}D_{0+}^{2H}h,\ H<1/2.
\eeqlb
\begin{eqnarray*}
(K_H^{-1}h)(s)=s^{H-1/2}I_{0+}^{1/2-H}s^{1/2-H}h' .
\end{eqnarray*}

In the present paper, we consider the following SDE driven by fractional Brownian motion with $H>1/2$:
\beqlb\label{2.0}
\d X_t=b(t,X_t)\d t+\sigma(t)\d B^H_t,\ X_0=x\in\R^d,\ t\in[0,T],
\eeqlb
where $b:[0,T]\times\R^d\rightarrow\R^d,\ \sigma:[0,T]\rightarrow\R^d\times\R^d$.\\
In \cite{Nualart&Rascanu02a}, the authors proved the existence and uniqueness result of the solution to \eqref{2.0} and showed the solution has finite moments.
Define $P_tf(x)=\E f(X_t^x), \ t\in[0,T], \ f\in\B_b(\R^d)$,
where $X_t^x$ is the solution to \eqref{2.0} with $X_0=x$ and $\B_b(\R^d)$ denotes the set of all bounded measurable functions on $\R^d$.
The purpose of this paper is to establish the integration by parts formula and present some applications.
To conclude this section, for $0<\alpha\leq1$, let $C^\lambda([0,T];\mathbb{R}^d)$ be the space of $\alpha$-H\"{o}lder continuous functions $f:[0,T]\rightarrow\mathbb{R}^d$
and set
$$\|f\|_\infty:=\sup_{0\leq t\leq T}|f(t)|,\ \ \|f\|_\alpha:=\sup\limits_{0\leq s<t\leq T}\frac{|f(t)-f(s)|}{|t-s|^\alpha}.$$

\section{Main result and its proof}

\setcounter{equation}{0}

To start with, let
$$\nabla_yf(x)=\lim\limits_{\epsilon\downarrow0}\frac{f(x+\epsilon y)-f(x)}{\epsilon}$$
and make the following assumptions on the coefficients: $\mathbf{(A)}$

\begin{itemize}
\item[(i)] there exist positive constants $L_i, i=1,2,3$ such that\\
a) H\"{o}lder continuous in time of order $H-1/2<\gamma\leq1$:
$$|b(t,x)-b(s,x)|+|\nabla b(t,\cdot)(x)-\nabla b(s,\cdot)(x)|\leq L_1|t-s|^\gamma,\ \forall t,s\in[0,T], \ x\in\R^d ;$$
b) H\"{o}lder continuous of order $1-1/(2H)<\rho<1$:
$$|\nabla b(t,\cdot)(x)-\nabla b(t,\cdot)(y)|\leq L_2|x-y|^\rho,\ \forall x,y\in\R^d,\ t\in[0,T];$$
c) boundedness:
$$|\nabla b(t,\cdot)(x)|\leq L_3,\ \forall x\in\R^d,\ t\in[0,T].$$

\par

\item[(ii)] $\sigma$ is H\"{o}lder continuous of order $(1-H)\vee(H-1/2)<\delta\leq1$ with positive constant $K$:
$$|\sigma(t)-\sigma(s)|\leq K|t-s|^{\delta},\ \forall t,s\in[0,T],$$
and $\sigma^{-1}$ is bounded.
\end{itemize}
According to \cite[Theorem 2.1]{Nualart&Rascanu02a}, the condition $\mathbf{(A)}$ ensures the equation \eqref{2.0} a unique solution,
whose sample paths are H\"{o}lder continuous of order less than $H$.\\
Set
\beqnn
\ar M_t\ar=\frac{1}{\Gamma(3/2-H)T}
\Bigg\{\int_0^t\left\langle s^{\frac{1}{2}-H}\sigma^{-1}(s)[y-s\nabla_{y}b(s,\cdot)(X_s)],\d W_s\right\rangle\cr
\ar\ar+\left(H-\frac{1}{2}\right)\int_0^t\left\langle s^{H-\frac{1}{2}}\sigma^{-1}(s)[y-s\nabla_{y}b(s,\cdot)(X_s)]
       \int_0^s\frac{s^{\frac{1}{2}-H}-r^{\frac{1}{2}-H}}{(s-r)^{\frac{1}{2}+H}}\d r,\d W_s\right\rangle\cr
\ar\ar-\left(H-\frac{1}{2}\right)\int_0^t\left\langle s^{H-\frac{1}{2}}\int_0^s
       \frac{s\sigma^{-1}(s)\nabla_{y}b(s,\cdot)(X_s)-r\sigma^{-1}(r)\nabla_{y}b(r,\cdot)(X_r)}{(s-r)^{\frac{1}{2}+H}}r^{\frac{1}{2}-H}\d r,\d W_s\right\rangle\cr
\ar\ar+\left(H-\frac{1}{2}\right)\int_0^t\left\langle s^{H-\frac{1}{2}}\int_0^s \frac{\sigma^{-1}(s)-\sigma^{-1}(r)}{(s-r)^{\frac{1}{2}+H}}r^{\frac{1}{2}-H}\d ry,\d W_s\right\rangle\Bigg\},\ \ \ t\in[0,T].
\eeqnn

The main result in this section is the following.

\btheorem\label{T3.1}
Assume $\mathbf{(A)}$ and let $y\in\R^d$ be fixed.
For any $f\in C_b^1(\R^d)$, there holds
\beqnn
P_T(\nabla_y f)(x)=\E\left[f(X_T^x)M_T\right].
\eeqnn
\etheorem

To prove the theorem, let $X_t$ solve the equation $\eqref{2.0}$ with the initial data $x$ and for any $\epsilon\in[0,1]$,
let $X_t^\epsilon$ solve the following equation
\beqlb\label{3.1a}
\d X_t^\epsilon=b(t,X_t)\d t+\sigma(t)\d B_t^H+\frac{1}{T}\epsilon y\d t,\ X_0^\epsilon=x\in\R^d, \ t\in[0,T].
\eeqlb
Then it is clear that $X_t^\epsilon=X_t+\frac{t}{T}\epsilon y,\ t\in[0,T]$. In particular, $X_T^\epsilon=X_T+\epsilon y$.\\
To see that $(X,X^\epsilon)$ is a coupling by change of measure for the solution to $\eqref{2.0}$,
we need to reformulate the equation $\eqref{3.1a}$ by using a new fractional Brownian motion.
To this end, let
\beqnn
\eta(t)=b(t,X_t)-b(t,X_t^\epsilon)+\frac{1}{T}\epsilon y,\ t\in[0,T]
\eeqnn
and
\beqnn
\tilde{B}^H_t&=&B^H_t+\int_0^t\sigma^{-1}(s)\eta(s)\d s\\
&=&\int_0^tK_H(t,s)\left(\d W_s+K_H^{-1}\left(\int_0^\cdot\sigma^{-1}(r)\eta(r)\d r\right)(s)\d s\right)\\
&=:&\int_0^tK_H(t,s)\d\tilde{W}_s.
\eeqnn
Now, let
\begin{eqnarray}\label{3.2a}
R_\epsilon&=&\exp\Big[-\int_0^TK_H^{-1}\left(\int_0^\cdot\sigma^{-1}(r)\eta(r)\d r\right)(t)\d W_t\cr
&&~~~~~~~-\frac{1}{2}\int_0^T\left|K_H^{-1}\left(\int_0^\cdot\sigma^{-1}(r)\eta(r)\d r\right)(t)\right|^2\d t\Big].
\end{eqnarray}
The next two results provide the exponential integrability of the r.v.
$\int_0^T|K_H^{-1}(\int_0^\cdot\sigma^{-1}(r)\eta(r)\d r)(t)|^2\d t$ shown in \eqref{3.2a}
and the convergence of $R_\epsilon$ in the $L^1(\P)$ sense, respectively.

\bproposition\label{P3.1}
Suppose that $\mathbf{(A)}$ holds.
Then, for any $\theta\in\R^+$, we have
\beqnn
\E\exp\left[\theta\int_0^T\left|K_H^{-1}\left(\int_0^\cdot\sigma^{-1}(r)\eta(r)\d r\right)(s)\right|^2\d s\right]<\infty.
\eeqnn
\eproposition

\bproposition\label{P3.2}
Suppose that $\mathbf{(A)}$ holds.
Then, there holds in $L^1(\P)$
\beqnn
\lim\limits_{\epsilon\rightarrow0}\frac{R_\epsilon-1}{\epsilon}=-M_T.
\eeqnn
\eproposition

To prove Proposition \ref{P3.1} and Proposition \ref{P3.2}, we state the following lemma concerning the estimations of the solution $X$.
The proof is modified from the one proposed in \cite{Hu&Nualart07} (see also \cite{Saussereau12}) and so we omit it.

\blemma\label{L3.1}
Assume $\mathbf{(A)}$.
Then, there hold
\begin{eqnarray*}
\|X\|_\infty&\leq&C\e^{L_3T}\left[1+T+T^{\gamma+1}+\|B^H\|_\lambda\left(T^\lambda+T^{\lambda+\delta}\right)\right]\\
&\leq&C(T)\left(1+\|B^H\|_\lambda\right),
\end{eqnarray*}
and for any $t,s\in[0,T]$,
\begin{eqnarray*}
|X_t-X_s|&\leq& C\Big\{\left[1+T^\gamma+\e^{L_3T}(1+T+T^{\gamma+1})\right]|t-s|+\e^{L_3T}\|B^H\|_\lambda\left(T^\lambda+T^{\lambda+\delta}\right)|t-s|\\
&&+\|B^H\|_\lambda\left(|t-s|^\lambda+|t-s|^{\lambda+\delta}\right)\Big\}\\
&\leq&C(T)\left(|t-s|+\|B^H\|_\lambda|t-s|^\lambda\right),
\end{eqnarray*}
where and in what follows, $C$ denotes a generic constant, $\lambda$ is chosen satisfying $1-\delta<\lambda<H$ and $\lambda\rho>H-1/2$.
\elemma

\textbf{Proof of Proposition \ref{P3.1}.}
We first note that $K_H^{-1}\left(\int_0^\cdot\sigma^{-1}(r)\eta(r)\d r\right)$ is well-defined.\\
Indeed, due to $\mathbf{(A)}$ and the relation $X_t^\epsilon=X_t+\frac{t}{T}\epsilon y,\ t\in[0,T]$, we get
\beqnn
|\eta(t)-\eta(s)|\leq C(|t-s|^\gamma+|t-s|+|X_t-X_s|).
\eeqnn
So, there holds $\eta(\cdot)\in C^{\min\{\gamma,H-\epsilon\}}([0,T],\R^d)$.
As a consequence, the fact that
\beqnn
|\sigma^{-1}(t)\eta(t)-\sigma^{-1}(s)\eta(s)|\leq C\left(|\eta(t)-\eta(s)|+|t-s|^\delta\right),
\eeqnn
and the condition $\gamma,\delta>H-1/2$ imply that $\sigma^{-1}\eta$ is H\"{o}lder continuous of order larger than $H-1/2$.
Then we have $\sigma^{-1}\eta\in I_{0+}^{H-1/2}(L^2([0,T],\R^d))$ and moreover, $\int_0^\cdot\sigma^{-1}(r)\eta(r)\d r\in I_{0+}^{H+1/2}(L^2([0,T],\R^d))$.

Next, we consider the exponential integrability of the r.v. $\int_0^T|K_H^{-1}(\int_0^\cdot\sigma^{-1}(r)\eta(r)\d r)(t)|^2\d t$.\\
By \eqref{2.0'}, we obtain
\beqlb\label{P3.1-1}
&&K_H^{-1}\left(\int_0^\cdot\sigma^{-1}(r)\eta(r)\d r\right)(s)
=s^{H-\frac{1}{2}}D_{0+}^{H-\frac{1}{2}}\left(r^{\frac{1}{2}-H}\sigma^{-1}(r)\eta(r)\right)(s)\cr
&=&\frac{1}{\Gamma(\frac{3}{2}-H)}\Bigg[s^{\frac{1}{2}-H}\sigma^{-1}(s)\eta(s)
   +\left(H-\frac{1}{2}\right)s^{H-\frac{1}{2}}\sigma^{-1}(s)\eta(s)\int_0^s\frac{s^{\frac{1}{2}-H}-r^{\frac{1}{2}-H}}{(s-r)^{\frac{1}{2}+H}}\d r\cr
&&~~~~~~~~~~~~~~~+\left(H-\frac{1}{2}\right)s^{H-\frac{1}{2}}\int_0^s\frac{\sigma^{-1}(s)\eta(s)-\sigma^{-1}(r)\eta(r)}{(s-r)^{\frac{1}{2}+H}}r^{\frac{1}{2}-H}\d r\Bigg]\cr
&=:&\frac{1}{\Gamma(\frac{3}{2}-H)}[I_1+I_2+I_3].
\eeqlb
Noting that
$$\int_0^s\frac{s^{\frac{1}{2}-H}-r^{\frac{1}{2}-H}}{(s-r)^{\frac{1}{2}+H}}\d r
=\int_0^1\frac{u^{\frac{1}{2}-H}-1}{(1-u)^{\frac{1}{2}+H}}\d u\cdot s^{1-2H}<\infty,$$
we easily get
\beqlb\label{P3.1-2}
|I_1|+|I_2|\leq C\epsilon s^{\frac{1}{2}-H}|y|.
\eeqlb
Now, we focus on the term $I_3$.\\
Firstly, observe that
\beqlb\label{P3.1-3}
&&\left|\int_0^s\frac{\sigma^{-1}(s)\eta(s)-\sigma^{-1}(r)\eta(r)}{(s-r)^{\frac{1}{2}+H}}r^{\frac{1}{2}-H}\d r\right|\cr
&\leq&|\eta(s)|\int_0^s\frac{|\sigma^{-1}(s)-\sigma^{-1}(r)|}{(s-r)^{\frac{1}{2}+H}}r^{\frac{1}{2}-H}\d r
   +C\int_0^s\frac{|\eta(s)-\eta(r)|}{(s-r)^{\frac{1}{2}+H}}r^{\frac{1}{2}-H}\d r\cr
&\leq&C\left(\epsilon s^{\delta-2H+1}|y|+\int_0^s\frac{|\eta(s)-\eta(r)|}{(s-r)^{\frac{1}{2}+H}}r^{\frac{1}{2}-H}\d r\right).
\eeqlb
For the integral of the last inequality in \eqref{P3.1-3},
in view of the fundamental theorem for Bochner integral, $\mathbf{(A)}$ and Lemma \ref{L3.1}, we arrive at
\beqlb\label{P3.1-4}
\int_0^s\frac{|\eta(s)-\eta(r)|}{(s-r)^{\frac{1}{2}+H}}r^{\frac{1}{2}-H}\d r &=&\int_0^s\frac{r^{\frac{1}{2}-H}}{(s-r)^{\frac{1}{2}+H}}\left|b(s,X_s)-b(s,X_s^\epsilon)-(b(r,X_r)-b(r,X_r^\epsilon))\right|\d r\cr
&\leq&\int_0^s\frac{r^{\frac{1}{2}-H}}{(s-r)^{\frac{1}{2}+H}}
\bigg|\int_0^1\nabla b(s,\cdot)(X_s+u(X_s^\epsilon-X_s))(X_s^\epsilon-X_s)\d u\cr
&&~~~~~~~~~~~~~~~~~~~~~-\int_0^1\nabla b(r,\cdot)(X_r+u(X_r^\epsilon-X_r))(X_r^\epsilon-X_r)\d u\bigg|\d r\cr
&\leq&C\epsilon|y|\int_0^s\frac{r^{\frac{1}{2}-H}}{(s-r)^{\frac{1}{2}+H}}
\bigg[\int_0^1\left|X_s-X_r+u\frac{r-s}{T}\epsilon y\right|^\rho\d u\cr
&&~~~~~~~~~~~~~~~~~~~~~~~~~~~~~+|s-r|^\gamma+\frac{s-r}{T}\bigg]\d r\cr
&\leq&C\epsilon|y|\bigg[s^{2-2H}+s^{\gamma-2H+1}+s^{\rho-2H+1}+s^{\rho-2H+1}|y|^\rho\cr
&&~~~~~~~~+\|B^H\|_\lambda^\rho\left(s^{\rho-2H+1}+s^{\lambda\rho-2H+1}+s^{(\lambda+\delta)\rho-2H+1}\right)\bigg].
\eeqlb
Therefore, combining \eqref{P3.1-3} with \eqref{P3.1-4} yields
\beqlb\label{P3.1-5}
|I_3|&\leq&C\epsilon|y|\bigg[s^{\frac{3}{2}-H}+s^{\delta-H+\frac{1}{2}}+s^{\gamma-H+\frac{1}{2}}+s^{\rho-H+\frac{1}{2}}+s^{\rho-H+\frac{1}{2}}|y|^\rho\cr
&&~~~~~~~~+\|B^H\|_\lambda^\rho\left(s^{\rho-H+\frac{1}{2}}+s^{\lambda\rho-H+\frac{1}{2}}+s^{(\lambda+\delta)\rho-H+\frac{1}{2}}\right)\bigg].
\eeqlb
Then substituting \eqref{P3.1-2} and \eqref{P3.1-5} into \eqref{P3.1-1}, we obtain
\begin{eqnarray*}
\left|K_H^{-1}\left(\int_0^\cdot\sigma^{-1}(r)\eta(r)\d r\right)(s)\right|
&\leq&C\epsilon|y|\bigg[s^{\frac{1}{2}-H}+s^{\frac{3}{2}-H}+s^{\delta-H+\frac{1}{2}}+s^{\gamma-H+\frac{1}{2}}\cr
&&~~~~~~~~+s^{\rho-H+\frac{1}{2}}+s^{\rho-H+\frac{1}{2}}|y|^\rho
+\|B^H\|_\lambda^\rho\cr
&&~~~~~~~~\times
\left(s^{\rho-H+\frac{1}{2}}+s^{\lambda\rho-H+\frac{1}{2}}+s^{(\lambda+\delta)\rho-H+\frac{1}{2}}\right)\bigg].
\end{eqnarray*}
So, we have, for each $\theta\in\R^+$,
\beqlb\label{P3.1-6}
\theta\int_0^T\left|K_H^{-1}\left(\int_0^\cdot\sigma^{-1}(r)\eta(r)\d r\right)(s)\right|^2\d s
\leq C\epsilon^2|y|^2\left(1+|y|^{2\rho}+\|B^H\|_\lambda^{2\rho}\right).
\eeqlb
As a consequence, the Fernique theorem implies the desired result.
\fin

\textbf{Proof of Proposition \ref{P3.2}.}
Let $R_\epsilon=\exp\left[M_T^\epsilon-\frac{1}{2}\langle M^\epsilon\rangle_T\right]$.
Without lost of generality, we suppose $\epsilon\leq 1$.
We first claim that
\begin{eqnarray}\label{P3.2-1}
\lim\limits_{\epsilon\rightarrow0}\E\frac{R_\epsilon-1}{\epsilon}
=\lim\limits_{\epsilon\rightarrow0}\E\frac{M_T^\epsilon-\frac{1}{2}\langle M^\epsilon\rangle_T}{\epsilon}.
\end{eqnarray}
Indeed, by the elementary inequalities: $|e^x-1-x|\leq x^2e^{|x|}, x^2\leq e^{|x|}, \forall x\in\mathbb{R}$, we obtain
\beqnn
\left|\frac{R_\epsilon-1-(M_T^\epsilon-\frac{1}{2}\langle M^\epsilon\rangle_T)}{\epsilon}\right|
\ar\leq\ar\frac{1}{\epsilon}\left(M_T^\epsilon-\frac{1}{2}\langle M^\epsilon\rangle_T\right)^2\exp\left[|M_T^\epsilon|+\frac{1}{2}\langle M^\epsilon\rangle_T\right]\cr
\ar=\ar\epsilon^\frac{1}{2}\left(\frac{1}{\epsilon^{3/4}}M_T^\epsilon-\frac{1}{2\epsilon^{3/4}}\langle
       M^\epsilon\rangle_T\right)^2\exp\left[|M_T^\epsilon|+\frac{1}{2}\langle M^\epsilon\rangle_T\right]\\
\ar\leq\ar\epsilon^\frac{1}{2}\exp\left[\frac{1}{\epsilon^{3/4}}|M_T^\epsilon|+\frac{1}{2\epsilon^{3/4}}\langle
       M^\epsilon\rangle_T+|M_T^\epsilon|+\frac{1}{2}\langle M^\epsilon\rangle_T\right]\\
\ar\leq\ar\epsilon^\frac{1}{2}\exp\left[\frac{2}{\epsilon^{3/4}}|M_T^\epsilon|+\frac{1}{\epsilon^{3/4}}\langle M^\epsilon\rangle_T\right]\cr
\ar\leq\ar\epsilon^\frac{1}{2}\left(\exp\left[\frac{2}{\epsilon^{3/4}}M_T^\epsilon\right]+\exp\left[\frac{-2}{\epsilon^{3/4}}M_T^\epsilon\right]\right)
  \exp\left[\frac{1}{\epsilon^{3/4}}\langle M^\epsilon\rangle_T\right]\cr
\ar=\ar\epsilon^\frac{1}{2}\left(\exp\left[\frac{2}{\epsilon^{3/4}}M_T^\epsilon-\frac{4}{\epsilon^{3/2}}\langle M^\epsilon\rangle_T\right]
  +\exp\left[\frac{-2}{\epsilon^{3/4}}M_T^\epsilon-\frac{4}{\epsilon^{3/2}}\langle M^\epsilon\rangle_T\right]\right)\cr
\ar\ar\times\exp\left[\frac{1}{\epsilon^{3/4}}\left(1+\frac{4}{\epsilon^{3/4}}\right)\langle M^\epsilon\rangle_T\right].
\eeqnn
This, together with the H\"{o}lder inequality, \eqref{P3.1-6} and the Fernique theorem, implies that for small enough $\epsilon$,
\beqnn
\ar\ar\E\left|\frac{R_\epsilon-1-(M_T^\epsilon-\frac{1}{2}\langle M^\epsilon\rangle_T)}{\epsilon}\right|\cr
\ar\leq\ar\epsilon^\frac{1}{2} \left\{2\E\exp\left[\frac{4}{\epsilon^{3/4}}M_T^\epsilon-\frac{8}{\epsilon^{3/2}}\langle M^\epsilon\rangle_T\right]
  +2\E\exp\left[\frac{-4}{\epsilon^{3/4}}M_T^\epsilon-\frac{8}{\epsilon^{3/2}}\langle M^\epsilon\rangle_T\right]\right\}^{\frac{1}{2}}\cr
\ar\ar\times\left\{\E\exp\left[\frac{2}{\epsilon^{3/4}}\left(1+\frac{4}{\epsilon^{3/4}}\right)\langle M^\epsilon\rangle_T\right]\right\}^{\frac{1}{2}}\cr
\ar=\ar(4\epsilon)^\frac{1}{2}
\left\{\E\exp\left[\frac{2}{\epsilon^{3/4}}\left(1+\frac{4}{\epsilon^{3/4}}\right)
\langle M^\epsilon\rangle_T\right]\right\}^\frac{1}{2},
\eeqnn
which shows that \eqref{P3.2-1} is true.

Now, by \eqref{P3.2-1} and \eqref{P3.1-6}, we have
\begin{eqnarray*}
\lim\limits_{\epsilon\rightarrow 0}\mathbb{E}\frac{R_\epsilon-1}{\epsilon}=\lim\limits_{\epsilon\rightarrow 0}\mathbb{E}\frac{M_T^\epsilon}{\epsilon}.
\end{eqnarray*}
Observe that by \eqref{P3.1-1}, $M_T^\epsilon$ can be writen as
\beqnn
M_T^\epsilon=-\frac{1}{\Gamma(\frac{3}{2}-H)}\int_0^T\langle I_1+I_2+I_3,\d W_s\rangle=:\frac{1}{\Gamma(\frac{3}{2}-H)}(J_1+J_2+J_3).
\eeqnn
For $J_1$, it follows from the B-D-G inequality that
\beqlb\label{P3.2-2}
&&\E\sup\limits_{0\leq t\leq T}\left|
   \int_0^t\left\langle s^{\frac{1}{2}-H}\sigma^{-1}(s)\left[\frac{-\eta(s)}{\epsilon}
   +\frac{y-s\nabla_{y}b(s,\cdot)(X_s)}{T}\right],\d W_s\right\rangle\right|\cr
&\leq&
\E\left[\int_0^Ts^{1-2H}|\sigma^{-1}(s)|^2\left|\frac{b(s,X_s^\epsilon)-b(s,X_s)-\nabla_y b(s,\cdot)(X_s)\frac{s}{T}\epsilon}{\epsilon}\right|^2\d s\right]^\frac{1}{2}\rightarrow 0,
\eeqlb
as $\epsilon$ goes to 0.\\
Similarly, for $J_2$ we conclude that, as $\epsilon$ tends to 0,
\beqlb\label{P3.2-3}
&&\E\sup\limits_{0\leq t\leq T}\left|\int_0^t
\left\langle s^{H-\frac{1}{2}}\sigma^{-1}(s)\left[\frac{-\eta(s)}{\epsilon}
+\frac{y-s\nabla_{y}b(s,\cdot)(X_s)}{T}\right]
\int_0^s\frac{s^{\frac{1}{2}-H}-r^{\frac{1}{2}-H}}{(s-r)^{\frac{1}{2}+H}}\d r,\d W_s\right\rangle\right|\cr
&&\rightarrow 0.
\eeqlb
For $J_3$, we first observe that
\beqlb\label{P3.2-4}
&&\E\sup\limits_{0\leq t\leq T}\Bigg|\int_0^t
\left\langle -s^{H-\frac{1}{2}}\int_0^s\frac{\sigma^{-1}(s)\eta(s)-\sigma^{-1}(r)\eta(r)}{\epsilon(s-r)^{\frac{1}{2}+H}}r^{\frac{1}{2}-H}\d r,\d W_s\right\rangle\cr
&&-\int_0^t\left\langle s^{H-\frac{1}{2}}\int_0^s
       \frac{s\sigma^{-1}(s)\nabla_{y}b(s,\cdot)(X_s)-r\sigma^{-1}(r)\nabla_{y}b(r,\cdot)(X_r)}{T(s-r)^{\frac{1}{2}+H}}r^{\frac{1}{2}-H}\d r,\d W_s\right\rangle\cr
&&+\int_0^t\left\langle s^{H-\frac{1}{2}}\int_0^s \frac{\sigma^{-1}(s)-\sigma^{-1}(r)}{T(s-r)^{\frac{1}{2}+H}}r^{\frac{1}{2}-H}\d ry,\d W_s\right\rangle
   \Bigg|\cr
&\leq&\E\Bigg[\int_0^Ts^{2H-1}\Bigg|\int_0^s
   \frac{\sigma^{-1}(s)-\sigma^{-1}(r)}{(s-r)^{\frac{1}{2}+H}}\frac{b(s,X_s^\epsilon)-b(s,X_s)-\frac{s\epsilon}{T}\nabla_y b(s,\cdot)(X_s)}{\epsilon} r^{\frac{1}{2}-H}\cr
&&~~~~~~~~~~~~~~~~~~~~~+\sigma^{-1}(r)\frac{b(s,X_s^\epsilon)-b(s,X_s)-(b(r,X_r^\epsilon)-b(r,X_r))}{\epsilon(s-r)^{\frac{1}{2}+H}}r^{\frac{1}{2}-H}\cr
&&~~~~~~~~~~~~~~~~~~~~~-\sigma^{-1}(r)\frac{s\nabla_y b(s,\cdot)(X_s)-r\nabla_y b(r,\cdot)(X_r)}{T(s-r)^{\frac{1}{2}+H}}r^{\frac{1}{2}-H}\d r
  \Bigg|^2\d s\Bigg]^\frac{1}{2}\cr
  &\leq&C\overline{\lim\limits_{\epsilon\rightarrow0}}\E\Bigg[\int_0^Ts^{2H-1}\Bigg|\int_0^s
  \sigma^{-1}(r)\frac{b(s,X_s^\epsilon)-b(s,X_s)-(b(r,X_r^\epsilon)-b(r,X_r))}{\epsilon(s-r)^{\frac{1}{2}+H}}r^{\frac{1}{2}-H}\cr
&&~~~~~~~~~~~~~~~~~~~~~-\sigma^{-1}(r)\frac{s\nabla_y b(s,\cdot)(X_s)-r\nabla_y b(r,\cdot)(X_r)}{T(s-r)^{\frac{1}{2}+H}}r^{\frac{1}{2}-H}\d r
  \Bigg|^2\d s\Bigg]^\frac{1}{2}\cr
&=:&C\overline{\lim\limits_{\epsilon\rightarrow0}}\E\left(\int_0^Ts^{2H-1}h(s)\d s\right)^\frac{1}{2}.
\eeqlb
By Lemma \ref{L3.1} and the argument of Proposition \ref{P3.1}, we have
\beqnn
\lim\limits_{\epsilon\rightarrow0}h(s)=0
\eeqnn
and
\beqnn
s^{2H-1}h(s)&\leq&C\bigg[s^{3-2H}+s^{2\gamma-2H+1}+s^{2\rho-2H+1}+s^{2\rho-2H+1}|y|^{2\rho}\cr
&&~~~+\|B^H\|_\lambda^{2\rho}\left(s^{2\rho-2H+1}+s^{2\lambda\rho-2H+1}+s^{2(\lambda+\delta)\rho-2H+1}\right)\bigg].
\eeqnn
Then the dominated convergence theorem implies
\begin{eqnarray}\label{P3.2-5}
\overline{\lim\limits_{\epsilon\rightarrow0}}\E\left(\int_0^Ts^{2H-1}h(s)\d s\right)^\frac{1}{2}=0.
\end{eqnarray}
So, by \eqref{P3.2-2}-\eqref{P3.2-5}, we complete the proof.
\fin

We now turn to the proof of Theorem \ref{T3.1} itself.

\textbf{Proof of Theorem \ref{T3.1}.}
Proposition \ref{P3.1} ensures that $\{\tilde{B}_t\}_{t\in[0,T]}$ is a $d$-dimensional fractional Brownian motion under the
probability $R_\epsilon\d\P$ by the Girsanov theorem for the fractional Brownian motion (see e.g., \cite[Theorem 4.9]{Decreusefond&Ustunel98a} or \cite[Theorem 2]{Nualart&Ouknine02b}).
Rewrite \eqref{3.1a} as follows
\beqnn
\d X_t^\epsilon=b(t,X_t^\epsilon)\d t+\sigma(t)\d\tilde{B}^H_t,\ X_0^\epsilon=x.
\eeqnn
Consequently, $(X,X^\epsilon)$ is a coupling by change of measure with changed probability $R_\epsilon\P$.
Moreover, since $R_0=1$, by \cite[Theorem 2.1]{Wang12b} and Proposition \ref{P3.2}, we derive the desired result.
\fin

\section{Some applications: shift Harnack type inequalities and absolute continuity of the law}

\setcounter{equation}{0}

In this section, we give some applications of Driver type integration by parts formula for $P_t$.

\btheorem\label{T4.1}
Assume $\mathbf{(A)}$ and let $y\in\R^d$ be fixed.
Then there exist constants $a(T)$ and $b(T)$ such that
\begin{itemize}
\item[(1)] for any nonnegative $f\in\mathscr{B}_b(\R^d)$,
 $$(P_Tf)^p\leq\left(P_T{\{f(y+\cdot)\}}^p\right)\exp\left\{C\frac{p}{p-1}
\left[1+a(T)+b(T)\left(1\vee\frac{p|y|}{p-1}\right)^{\frac{2\rho}{1-\rho}}\right]|y|^2\right\}.$$

\par

\item[(2)] for any positive $f\in\mathscr{B}_b(\R^d)$,
$$P_T\log f\leq\log P_T\{f(y+\cdot)\}
+C\left[1+a(T)+b(T)\left(1\vee\frac{p|y|}{p-1}\right)^{\frac{2\rho}{1-\rho}}\right]|y|^2.$$

\end{itemize}
\etheorem

\textbf{Proof.}
By Theorem \ref{T3.1} and the Young inequality (see, for instance, \cite[Lemma 2.4]{Arnaudon&Thalmaier&Wang09a}),
we deduce that, for all $\theta>0$,
\begin{eqnarray}\label{T4.1-1}
|P_T(\nabla_y f)|-\theta\left[P_T(f\log f)-(P_T f)(\log P_T f)\right]
\ar\leq\ar\theta\log\mathbb{E}\exp\left[\frac{1}{\theta}M_T\right]\cdot P_Tf\cr
\ar\leq\ar\frac{\theta}{2}\log\mathbb{E}\exp\left[\frac{2}{\theta^2}\langle M\rangle_T\right]\cdot P_T f.
\end{eqnarray}
On the other hand, in view of the expression of $M_t$ and Lemma \ref{L3.1}, we conclude that
\beqnn
\langle M\rangle_T\leq C\left(a(T)+\tilde{a}(T)\|B^H\|^{2\rho}_\lambda\right)|y|^2,
\eeqnn
where
\beqnn
a(T)=\left\{1+T+T^2+T^{2\delta}+T^{2(\delta+1)}+T^{2(\gamma+1)}+[1+T^\gamma+\e^{L_3T}\left(|x|+T+T^{1+\gamma}\right)]^{2\rho}T^{2(\rho+1)}\right\}T^{-2H}
\eeqnn
and
\beqnn
\tilde{a}(T)=\left[1+T^{2\delta\rho}+\e^{2L_3\rho T}\left(1+T^{\delta}\right)^{2\rho}T^{2\rho}\right]T^{2(\lambda\rho-H+1)}.
\eeqnn
Then a similar argument to that in \cite[Lemma 3.7]{Fan12} shows that
\begin{eqnarray*}
\mathbb{E}\exp\left[\frac{2}{\theta^2}\langle M\rangle_T\right]
\ar\leq\ar\exp\left\{\frac{|y|^2}{\theta^2}C\left[1+a(T)+b(T)\left(1\vee\frac{p|y|}{p-1}\right)^{\frac{2\rho}{1-\rho}}\right]\right\},
\end{eqnarray*}
where $b(T)=\tilde{a}(T)^{\frac{1}{1-\rho}}$.\\
This, together with \eqref{T4.1-1}, yields
\begin{eqnarray*}
\ar\ar|P_T(\nabla_y f)|-\theta\left[P_T(f\log f)-(P_T f)(\log P_T f)\right]\cr
\ar\leq\ar
C\left[1+a(T)+b(T)\left(1\vee\frac{p|y|}{p-1}\right)^{\frac{2\rho}{1-\rho}}\right] \frac{|y|^2}{\theta}P_T f.
\end{eqnarray*}
Therefore, due to \cite[Proposition 2.3]{Wang12b}, it is easy to follow the desired result.
\fin

These inequalities above allow us to study the existence of distribution density of the solution. That is,

\bcorollary\label{C4.1}
Assume $\mathbf{(A)}$.
Then, for any $t>0$, the law of the solution $X_t$ of $\eqref{2.0}$ is absolutely continuous with respect to the Lebesgue measure.
\ecorollary

\textbf{Proof.}
Without lost of generality, we only consider the case $t=T$.
Let
$$h(y)=C\frac{p}{p-1}\left[1+a(T)+b(T)\left(1\vee\frac{p|y|}{p-1}\right)^{\frac{2\rho}{1-\rho}}\right]|y|^2.$$
By Theorem \ref{T4.1}, we deduce that, for any nonnegative $f\in\mathscr{B}_b(\R^d)$,
\begin{eqnarray}\label{C4.1-1}
(P_T f(x))^p{\e}^{-h(y)}\leq\left(P_T{\{f(y+\cdot)\}}^p\right)(x).
\end{eqnarray}
For any Lebesgue-null set $A\in\R^d$, choosing $f=I_A$ and integrating both sides with respect to $\d y$ in \eqref{C4.1-1} yield
$$(P_T\I_A(x))^p\int_{\R^d}e^{-h(y)}dy\leq\int_{\R^d}\int_{\R^d}\I_A(y+z)\d y\P\circ(X_T^x)^{-1}(\d z)=0.$$
Consequently, we have
$P_T\I_A(x)=0$, i.e., $\P\circ(X_T^x)^{-1}(A)=0$.
Then the proof is complete.
\fin

\end{document}